\documentclass[11pt]{article}
\usepackage{amssymb,amsfonts,amsmath,amsthm}
\usepackage{color} 
\usepackage[normalem]{ulem} 
\newtheorem{theorem}{Theorem}[section]
\newtheorem{thm}[theorem]{Theorem}
\newtheorem{prop}[theorem]{Proposition}
\newtheorem{lem}[theorem]{Lemma}
\newtheorem{rem}[theorem]{Remark}
\newtheorem{cor}[theorem]{Corollary}
\newtheorem{conj}[theorem]{Conjecture}

\makeatletter
   \@addtoreset{equation}{section}
\makeatother

\def\CT{\mathop{\mathrm{CT}}}
\def\ZZ{\mathbb{Z}}
\def\NN{\mathbb{N}}
\def\CC{\mathbb{C}}
\def\T{\overline{T}}
\begin{document}

\title{A Unified Elementary Approach to the Dyson, Morris, Aomoto, and Forrester Constant Term Identities}
\author{
  {\vspace{0.2cm}}
  {\small Ira M. Gessel$^1$, Lun Lv$^2$, Guoce Xin$^3$, Yue Zhou$^4$}\\
  {\small $^1$ Department of Mathematics}\\
  {\vspace{0.2cm}}
  {\small Brandeis University, Waltham, MA 02454-9110, USA}\\
  {\small $^{2,3,4}$Center for Combinatorics, LPMC-TJKLC}\\
  {\vspace{0.2cm}}
  {\small Nankai University, Tianjin 300071, P.R. China}\\
  { \small $^1$gessel@brandeis.edu\ \ \  \  $^2$lvlun@mail.nankai.edu.cn}\\
  { \small \ \ \ $^3$gxin@nankai.edu.cn\ \ \ $^4$zhouyue@mail.nankai.edu.cn}\\
 }
\date{February 21, 2008}
\maketitle

\begin{abstract}
We introduce an elementary method to give unified proofs of the
Dyson, Morris, and Aomoto identities for constant terms of Laurent
polynomials. These identities can be expressed as equalities of
polynomials and thus can be proved by verifying them for
sufficiently many values, usually at negative integers where they
vanish. Our method also proves some special cases of the Forrester
conjecture.\end{abstract}



\section{Introduction}

In 1962, Freeman Dyson \cite{dyson} conjectured the following
identity:

\begin{theorem}\label{t-dyson}
For nonnegative integers $a_0,a_1,\ldots ,a_n$,
\begin{equation}
\CT_{x} \prod_{0\le i\ne j \le n}
\left(1-\frac{x_i}{x_j}\right)^{\!\!a_j} =
 \frac{(a_0+a_1+\cdots+a_n)!}{a_0!\, a_1!\, \cdots a_n!},\label{e-dyson}
\end{equation}
where $\CT_x$ denotes the constant term.
\end{theorem}
Dyson's conjecture was proved independently  before his  paper was
published by Gunson~\cite{gunson} and  by Wilson \cite{wilson}, and
an elegant recursive proof was later found by Good \cite{good1}.

Similar identities for constant terms of Laurent polynomials expressed as products are of considerable interest, and we shall discuss several of them in this paper.

First is an identity of Morris \cite{morris82}.
For  $a, b, k\in \mathbb{N}$ (the nonnegative integers) define
\begin{equation}\label{hn}
H(x_0,x_1,\ldots,x_n;a,b,k) := \prod_{l=1}^n
\left(1-\frac{x_l}{x_0}\right)^{\!\!a}
\left(1-\frac{x_0}{x_l}\right)^{\!\!b}\prod_{1\leq i\neq j\leq n}
  \left( 1-\frac{x_i}{x_j} \right)^{\!\!k}.
\end{equation}

Morris proved the following result:

\begin{thm}\label{t-morris}
\begin{equation}\label{e-morris}
\CT_x H(x_0,x_1,\ldots,x_n;a,b,k)=M_n(a,b,k),
\end{equation}
where
\begin{equation}\label{mn}
M_n(a,b,k):=\prod_{l=0}^{n-1}\frac{(a+b+kl)!\,(k(l+1))!}{(a+kl)!\,(b+kl)!\,k!}.
\end{equation}
\end{thm}
Since $H(x_0,\dots,x_n;a,b,k)$ is homogeneous of degree $0$ in the
$x_i$, setting $x_0=1$ in \eqref{e-morris}  gives an equivalent
result, which is the form stated by Morris \cite{morris82}.

A generalization of the Morris identity was given by Aomoto
\cite{Aomoto87}, who extended Selberg's integral to obtain a formula
equivalent to the following constant term identity \cite{kadell98}:
Let
$$
A_m(x_0,x_1,\ldots,x_n;a,b,k):=\prod_{l=1}^{n}\Big(1-\frac{x_l}{x_0}\Big)^{\chi(l\leq
m)}H(x_0,x_1,\ldots,x_n;a,b,k),$$ where $\chi(S)$ equals $1$ if the
statement $S$ is true and $0$ otherwise.

\begin{thm}\label{t-aomoto}
\begin{align}\label{extmorris}
\qquad \CT_x A_m(x_0,x_1,\ldots,x_n;a,b,k)
=\prod_{l=0}^{n-1}\frac{[a+b+kl+\chi(l\geq
n-m)]!\,(k(l+1))!}{[a+kl+\chi(l\geq n-m)]!\,(b+kl)!\,k!}.
\end{align}
\end{thm}

Another generalization was conjectured by Forrester \cite{forr95}:
\begin{conj}\label{conjecture}
We have
\begin{multline}\label{e-forrester}
\CT_x\prod_{i,j=n_0+1,i\neq j}^{n}
   \left(1-\frac{x_i}{x_j}\right)H(x_0,x_1,\ldots,x_n;a,b,k)\\
=M_{n_0}(a,b,k)\prod_{j=0}^{n_1-1}\frac{(j+1)((k+1)j+a+b+kn_0)!\,(kj+k+j+kn_0)!}
{k!\,((k+1)j+a+kn_0)!\,((k+1)j+b+kn_0)!},
\end{multline}
where $n=n_0+n_1.$
\end{conj}

In \cite{forr95}, Forrester proved the special case $a=b=0$ (for all
$k$, $n_0$, and $n_1$) using a formula due to Bressoud and Goulden
\cite{bressoud-goulden}, and the case $k=1$ (for all $a$, $b$,
$n_0$, and $n_1$). Kaneko \cite{kaneko2000, kaneko2002} proved the
special cases $n_1=2,3$ and $n_1=n-1$. Moreover, Forrester and Baker
\cite{baker-forr} formulated a $q$-analog of Conjecture
\ref{conjecture}, which was recently studied by Kaneko
\cite{kaneko00}.

Our objective in this paper is to introduce an elementary method
which leads to new proofs of the Dyson, Morris, and Aomoto
identities. Moreover, our method can be used to obtain some partial
results on Forrester's conjecture.

The idea behind the proofs is the well-known fact that to prove the
equality of two polynomials of degree at most $d$, it is sufficient
to prove that they are equal at $d+1$ points. This approach was used
by Dyson \cite{dyson} to prove the case $n=3$ of \eqref{e-dyson}.
Dyson used Dougall's method \cite{dougall}, in which most of the
points are obtained by induction, making heavy use of the symmetry
of \eqref{e-dyson}. This approach does not seem to generalize beyond
$n =4$. In our approach we use the fact that, as a polynomial in
$a_0$, the right side of \eqref{e-dyson} vanishes for $a_0 = -1, -2,
\dots, -(a_1+\cdots +a_n)$, and we show that the same is true of the
left side.

The same idea was used by Gessel and Xin
\cite{gessel-xin} in proving a $q$-analog of Theorem \ref{t-dyson},
which was conjectured by George Andrews \cite{andrews-qdyson} in
1975, and first proved by Zeilberger and Bressoud
\cite{zeil-bressoud} in 1985.

In all of the proofs, it is routine to show that after fixing all
but one parameter, the constant term is a polynomial of degree at
most $d$ in the remaining parameter, say $a$, and that the left side
agrees with the right side when $a=0$. The proofs then differ in
showing that both sides vanish at $d$ additional points. In
\eqref{e-morris}, \eqref{extmorris} and \eqref{e-forrester}, these
polynomials may have multiple roots. We use the polynomial approach
to prove the cases in which the roots are distinct, and use another
argument, based on the form of the constant term as a function of
all the parameters (Proposition \ref{p-rat}), to extend the result
to the general case.

\section{Polynomials, vanishing coefficients, and a rationality result}
\label{s2}
In this section we prove several lemmas that will be needed in the
proofs of the constant term identities. First, we show in Lemma
\ref{s-poly} that the constant terms in these identities can be
expressed as polyomials. Next, Lemma \ref{l-ineq} is useful in
showing that these polynomials vanish at certain negative integers.
Lemma \ref{l-dys-root}, which applies Lemma \ref{l-ineq} to
coefficients of the Dyson product, gives Dyson's conjecture and is
also needed in the proof of Proposition \ref{p-rat}, which allows us
to deal with polynomials with multiple roots.

\subsection{A polynomial characterization}
\label{s-poly} Fundamental to our approach is the following lemma,
which shows that the constant terms we study are polynomials.

\begin{lem}\label{lemp} Let $a_0,\dots, a_n$ be nonnegative integers, $d:=a_{1}+\cdots+a_{n}$,
and let $L(x_1,\dots,x_n)$ be a Laurent
polynomial independent of $a_{0}$. Then for fixed
$a_{1},\ldots,a_{n}$, the constant term
\begin{equation}
\label{e-Q}
Q(a_{0},a_{1},\ldots,a_{n}):=\CT_{x}x_0^{k_0}L(x_1,\dots,x_n)\prod_{l=1}^n
\left(1-\frac{x_l}{x_0}\right)^{\!\!a_{0}}
\left(1-\frac{x_0}{x_l}\right)^{\!\!a_{l}}
\end{equation}
is a polynomial in $a_{0}$ of degree at most $d+k_0$ for any integer
$k_0\ge -d$.
\end{lem}

\begin{proof}
We can rewrite $Q(a_{0},a_{1},\ldots,a_{n})$ in the following form
\begin{align*}
Q(a_{0},a_{1},\ldots,a_{n})=(-1)^{\sum_l a_{l}} \CT_{x}
\frac{\prod_{l=1}^{n}(x_{0}-x_{l})^{a_{0}+a_{l}}}{x_{0}^{na_{0}-k_0}\prod_{l=1}^{n}x_{l}^{a_{l}}}
L(x_{1},\ldots,x_{n}).
\end{align*}
Expanding each $(x_{0}-x_{l})^{a_{0}+a_{l}}$ as
$\sum_{i_{l}=0}^{\infty}(-1)^{i_l}{a_{0}+a_{l} \choose
i_{l}}x_{0}^{a_{0}+a_{l}-i_{l}}x_{l}^{i_{l}}$, we get
\begin{align}
Q(a_{0},& a_{1},\ldots,a_{n})\nonumber\\
 &=\CT_{x} \sum_{i_{1},\ldots,i_{n}}(-1)^{\sum_l a_{l}+i_{l}}
   {a_{0}+a_{1}\choose i_{1}}\cdots{a_{0}+a_{n}
    \choose i_{n}}x_{0}^{k_0+\sum (a_{l}-i_{l})}
    \prod_{l=1}^{n}x_{l}^{i_{l}-a_l}L(x_{1},\ldots,x_{n})\nonumber\\
 &=  \sum_{i_{1},\ldots,i_{n}}
 (-1)^{k_0}
   {a_{0}+a_{1}\choose i_{1}}\cdots{a_{0}+a_{n}\choose i_{n}}\cdot
   \CT_{x_{1},\ldots,x_{n}}\prod_{l=1}^{n}x_{l}^{i_{l}-a_{l}}
   L(x_{1},\ldots,x_{n}),\label{expand1}
\end{align}
where the sum ranges over all nonnegative integers $i_1,\dots,i_n$
such that $i_1+\cdots +i_n=d+k_0$. To show that the degree of
$Q(a_0,a_1,\dots,a_n)$ in $a_0$ is at most $d+k_0$, it suffices to
show that every term has degree in $a_0$ at most $d+k_0$. This
follows from the fact that $\prod_{l=1}^{n}x_{l}^{i_{l}-a_{l}}
L(x_{1},\ldots,x_{n})$ is a Laurent polynomial independent of
$a_{0}$, and the fact that the degree of ${a_{0}+a_{1}\choose
i_{1}}\cdots{a_{0}+a_{n}\choose i_{n}}$ in $a_{0}$ is $d+k_0$, since
${a_{0}+a_{l}\choose i_{l}}$ is a polynomial in $a_{0}$ of
degree~$i_{l}$.
\end{proof}

Some corollaries of Lemma \ref{lemp} are given in the Appendix (Section \ref{appendix}).

Since $Q(a_0,a_1,\dots, a_n)$ as defined in \eqref{e-Q} is a
polynomial in $a_0$, we can extend it to all integers $a_0$, not
just nonnegative integers. It is useful to extend the meaning of the
right side of \eqref{e-Q} so that \eqref{e-Q} holds for negative
integers $a_0$. Since $(1-x_l/x_0)^{a_0}$ for $l\ge1$ is not a
Laurent polynomial unless $a_0$ is a nonnegative integer, we must
expand it as a Laurent series, but since
\begin{equation*}
\left(1-\frac{x_l}{x_0}\right)^{a_0}\!\! =
\left(-\frac{x_l}{x_0}\right)^{\!\!a_0}\!\!
\left(1-\frac{x_0}{x_l}\right)^{a_0}
\end{equation*}
we might conceivably expand this expression either in powers of $x_l/x_0$ or of $x_0/x_l$.
To make the expansion well-defined, we need to specify the ring in which we work.
We recall that for a ring $R$, the ring $R((x_1, x_2,
\dots, x_n))$  of formal Laurent series in $x_1,\dots, x_n$ with
coefficients in $R$ is the set of all formal series in these
variables in which only finitely many negative powers of $x_j$
appear for each $j$. Then it is sufficient to work in the ring
$\CC((x_0))((x_1,\dots, x_n))$ of formal Laurent series in
$x_1,\dots, x_n$ with coefficients in $\CC((x_0))$.
Informally, we
may think of $x_0$ as larger than all the other variables, so that
$x_l/x_0$ is small for $l\ge 1$. Thus we have in this ring the expansion
\begin{equation*}
(x_0 - x_l)^{a_0}
  = x_0^{a_0}\left(1-\frac{x_l}{x_0}\right)^{a_0}\!\!
  =\sum_{i=0}^\infty \binom{a_0}{i} x_0^{a_0-i}(-x_l)^i
\end{equation*}
for all integers $a_0$;
the alternative expansion
\begin{equation*}
(x_0 - x_l)^{a_0}
  = (-1)^{a_0}(x_l-x_0)^{a_0}
  =(-1)^{a_0}\sum_{i=0}^\infty \binom{a_0}{i} x_l^{a_0-i}(-x_0)^i
\end{equation*}
is not valid in this ring unless $a_0$ is a nonnegative integer.

\subsection{Vanishing coefficients}

Our goal is to evaluate special cases of the polynomial
$Q(a_0,a_1,\dots, a_n)$ given by Lemma \ref{lemp} by finding some of
their zeroes. The following lemma helps us to accomplish this.

\begin{lem}
\label{l-ineq}
Let $u_{ij}$ for $1\le i < j \le n$ be nonnegative integers and let $m_l$ and $v_l$ for $1\le l\le n$ be integers. If the coefficient of $x_1^{m_1}x_2^{m_2}\cdots x_n^{m_n}$ in
\begin{equation}
\label{e-ineq}
\frac{\prod_{1\le i<j\le n}(x_i-x_j)^{u_{ij}}}{\prod_{l=1}^n (1-x_l)^{v_l}}
\end{equation}
is nonzero then for some subset $T$ of $[n] := \{1,2,\dots, n\}$ we
have
\begin{align}
\sum_{\substack{i,j\in T\\ i<j}}u_{ij}&\le \sum_{i\in T}v_i -|T|\label{e-2.8}\\[-20pt]
\intertext{and}\notag\\[-20pt]
\sum_{1\le i<j\le n}\!\!u_{ij}&\ge \sum_{\substack{i,j\in T\\ i<j}}u_{ij} +\sum_{i\in \T}(m_i+v_i),
\label{e-2.9}
\end{align}
where $\T=[n]\setminus T$.
\end{lem}

\begin{proof}
Applying the formula
\begin{equation*} \label{e-xij}
\left({x_i}-{x_j}\right)^{u_{ij}}=\bigl((1-x_j)-(1-x_i)\bigr)^{u_{ij}}=\sum_{\alpha_{ij}+\alpha_{ji}=u_{ij}}
(-1)^{\alpha_{ij}}\binom{u_{ij}}{\alpha_{ij}}(1-x_i)^{\alpha_{ij}}(1-x_j)^{\alpha_{ji}}
\end{equation*}
and expanding, we can write \eqref{e-ineq} as a linear combination of
terms of the form
\[
\prod_{l=1}^{n}(1-x_l)^{-v_l}\!\! \prod_{1\leq i<j\leq n}
(1-x_i)^{\alpha_{ij}}(1-x_j)^{\alpha_{ji}}
  = \prod_{l=1}^{n}(1-x_l)^{-v_l} \!\!\prod_{1\leq i\ne j\leq n}
(1-x_i)^{\alpha_{ij}}.
\]
Now let $\alpha_i := \sum_{j=1}^n \alpha_{ij}$, where
$\alpha_{ii}=0.$ Then we may write this product as
\begin{equation}
\label{e-inprod}
\prod_{i=1}^n ( 1-x_i)^{\alpha_i - v_i}.
\end{equation}
If the coefficient of
 $x_1^{m_1}\cdots x_n^{m_n}$ in \eqref{e-inprod} is nonzero, then for each $i$, either $\alpha_i -v_i<0$ or
$m_i\le \alpha_i - v_i$. So if the coefficient of $x_1^{m_1}\cdots
x_n^{m_n}$ in \eqref{e-ineq} is nonzero, there exist nonnegative
integers $\alpha_{ij}$ with $\alpha_{ij}+\alpha_{ji}=u_{ij}$ for
$i\ne j$ and $\alpha_{ii}=0$ and a subset $T$ of $[n]$ such that
\begin{alignat}{2}
\alpha_{i}&\le v_i-1,&\quad&\text{for $i\in T$},  \label{dyson2}\\
\alpha_{i}&\geq m_i+v_i,&\quad&\text{for $i\in
\T$},\label{dyson3}
\end{alignat}
where $\alpha_i = \sum_{j=1}^n \alpha_{ij}$ and $\T=[n]\setminus T$.

Then
\begin{equation}\sum_{i\in T}\alpha_{i}
   =\sum_{\substack{i\in T\\j\in [n]}}\alpha_{ij}
   \ge\sum_{i,j\in T}\alpha_{ij}
   =\sum_{\substack{i,j\in T\\i<j}}(\alpha_{ij}+\alpha_{ji})
   = \sum_{\substack{i,j\in T\\i<j}}u_{ij}.
  \label{dyson5}
\end{equation}
Similarly,
\begin{equation}\label{dyson1}
\alpha_1+\cdots+\alpha_n=\sum_{1\le i<j\le n}\!\! u_{ij}.
\end{equation}

Summing \eqref{dyson2} for $i\in T$ gives
\[\sum_{i\in T}\alpha_{i} \le \sum_{i\in T}v_{i} -|T|,\]
so by \eqref{dyson5},
\[\sum_{\substack{i,j\in T\\i<j}}u_{ij} \le  \sum_{i\in T}v_i-|T|,\]
which is \eqref{e-2.8}.

Summing \eqref{dyson3} for $i\in \T$ gives
\[\sum_{i\in \T}\alpha_{i}\ge \sum_{i\in \T}(m_i+v_i).\]
Thus by \eqref{dyson1} and \eqref{dyson5} we have
\begin{equation*}
\sum_{1\le i<j\le n}\!\!u_{ij}
   = \sum_{i\in T}\alpha_{i} + \sum_{i\in \T}\alpha_{i}
   \ge \sum_{\substack{i,j\in T\\i<j}}u_{ij} +\sum_{i\in \T}(m_i+v_i),
\end{equation*}
and this is \eqref{e-2.9}.
\end{proof}

Next, we apply Lemma \ref{l-ineq} to prove the vanishing of some coefficients related to the Dyson product.

\begin{lem}
\label{l-dys-root} Let $a_1,\dots, a_n$  be nonnegative integers and
let $d:=a_1+\cdots+a_n$. Let $k_0, k_1,\dots, k_n$ be integers and
let $k$ be the sum of the positive integers among $k_1,\dots,k_n$.
For a subset $T\subseteq[n]$ we define $\sigma(T):=\sum_{i\in
T}a_i$, and we set
\[J:=
\bigcup_{T\subset
[n]}\{\sigma(T)+1,\sigma(T)+2,\ldots,\sigma(T)+k\},\] where the
union is over proper subsets $T$ of $[n]$. Then for every $a_0$ with
$-a_0\in [d]\setminus J$, we have
\begin{align}\label{dy-6}
\CT_x x_0^{k_0} x_1^{k_1}\cdots x_n^{k_n} \prod_{0\le i\ne j \le n}
\left(1-\frac{x_i}{x_j}\right)^{\!\!a_j} = 0.
\end{align}
\end{lem}
\begin{proof}
First we note that this coefficient is well defined for any negative
integer $a_0$, as explained at the end of Section \ref{s-poly}.
Next, since the product in \eqref{dy-6} is homogeneous of degree 0
in $x_0, x_1,\dots, x_n$, the constant term does not change if we
set $x_0$ equal to 1, as long as $k_0+\cdots+k_n=0$. (Otherwise the
constant term is 0.) Setting $a_0=-h$ and simplifying, we need to
show that
\begin{equation}\label{dyson6}
\CT_{{x}} \frac{\prod_{1\leq i<j\leq n}(x_i-x_j)^{a_i+a_j}}
{x_1^{na_1-k_1}\cdots
x_n^{na_n-k_n}\prod_{l=1}^{n}(1-x_l)^{h-a_l}}=0, \quad  \text{for
$h\in [d]\setminus J$}.
\end{equation}
We prove the contrapositive: Suppose that $h\in [d]$ but the left side of \eqref{dyson6}  is not 0. We shall show that
$h\in J$; i.e., $\sigma(T) < h\le \sigma(T) +k$ for some $T\subset [n]$.

We apply Lemma \ref{l-ineq} with $u_{ij}=a_{i}+a_{j}$, $m_i=na_i - k_i$, and
$v_i = h - a_i$. Then  for some subset $T\subseteq[n]$
we have
\begin{align}
\sum_{\substack{i,j\in T\\ i<j}}(a_i+a_j)&\le \sum_{i\in T}(h-a_i)-|T| \label{e-2.8d}\\[-20pt]
\intertext{and}\notag\\[-20pt]
\sum_{1\le i<j\le n}\!\!(a_i+a_j)&\ge \sum_{\substack{i,j\in T\\ i<j}}(a_i+a_j)
 +\sum_{i\in \T}\bigl((n-1)a_i+h-k_i\bigr).
\label{e-2.9d}
\end{align}
Let $t=|T|$. Then \eqref{e-2.8d} may be written as
\begin{equation*}
(t-1)\sum_{i\in T}a_i\le th-t - \sum_{i\in T}a_i,
\end{equation*}
and this implies that for $T\ne\varnothing$,
\begin{equation}
\label{dyson4} \sum_{i\in T}a_i < h.
\end{equation}
But \eqref{dyson4} also holds for $T=\varnothing$, since $h\ge 1$. We note that by
\eqref{dyson4}, $T\ne[n]$, since $h\le d$.

Similarly, \eqref{e-2.9d} gives
\begin{equation*}
(n-1)\sum_{i=1}^n a_i \ge (t-1) \sum_{i\in T}a_i + \sum_{i\in \T}(n-1)a_i +(n-t)h -\sum_{i\in \T}k_i.
\end{equation*}
Taking all the terms in the $a_i$ to the left side gives
\[(n-t)\sum_{i\in T}a_{i}\geq (n-t)h-\sum_{i\in \T}k_{i}\]
so since $T\ne [n]$,
\begin{equation}
\label{dyson5a} h\le \sum_{i\in T}a_{i}+k.
\end{equation}
Thus by \eqref{dyson4} and \eqref{dyson5a},
\[  \sum_{i\in T}a_{i}<h\le \sum_{i\in T}a_{i}+k,\]
which completes the proof.
\end{proof}

Dyson's conjecture is an easy consequence of Lemma \ref{l-dys-root}:
\begin{proof}[Proof of Theorem {\rm\ref{t-dyson}}]

Fix $a_1,\dots,a_n \in \mathbb{N}$.
Denote by $D_L(a_0)$ and $D_R(a_0)$ the left and right sides of
\eqref{e-dyson}. It is routine to check that
\begin{itemize}
    \item[1.] Both $D_L(a_0)$ and $D_R(a_0)$ are polynomials in $a_0$ of degree at most
    $d$ (by Lemma \ref{lemp});
    \item[2.] $D_L(0)=D_R(0)$ (by induction on $n$);
    \item[3.] $D_R(a_0)$ vanishes when $a_0=-1,-2,\dots, -d$.
\end{itemize}
Now apply Lemma \ref{l-dys-root} with $k_0=k_1=\cdots
=k_n=0$, so $k=0$ and  $J=\varnothing$. Then $D_L(a_0)$
also vanishes when $a_0=-1,-2,\dots, -d$. The theorem then follows
since two polynomials of degree at most $d$ are equal if they agree
at $d+1$ distinct points.
\end{proof}

\subsection{A rationality result}

We denote by $D_n(x;a_0, a_1,\dots, a_n)$ the Dyson product
\begin{equation*}
\prod_{0\le i\ne j \le n}
\left(1-\frac{x_i}{x_j}\right)^{\!\!a_j}.
\end{equation*}

Good \cite{good1} used Lagrange interpolation to derive the
following recursion in his proof of the Dyson conjecture: for
$a_0,a_1,\dots,a_n\ge 1$, we have
\begin{equation}
\label{e-good-recursion} D_n(x;a_0,a_1,\dots,a_n)=\sum_{i=0}^n
D_n(x;a_0,\dots,a_{i-1},a_i-1,a_{i+1},\dots,a_n).
\end{equation}

Using this recursion, Sills and Zeilberger \cite{sills-zeilberger},
Sills \cite{sills}, and Lv et.~al \cite{LvXinZhou} found explicit
formulas for some of the other coefficients of the Dyson product.
Their results suggest the following proposition, which we will need
in our approach to the Morris,  Aomoto, and Forrester constant
terms.

\begin{prop}\label{p-rat}
For any Laurent polynomial $L(x_0,\dots,x_n)$ independent of the $a_i$,
\begin{equation}\label{e-rationality}
\CT_x L(x_0,\dots,x_n) D_n(x;a_0,\dots,a_n)
=R(a_0,\dots,a_n)\frac{(a_0+a_1+\cdots +a_n)!}{a_0!\,a_1!\cdots
a_n!}
\end{equation}
for some rational function $R(a_0,\dots,a_n)$ of $a_0,\dots,a_n$.
\end{prop}


\begin{proof}
We proceed by induction on $n$. The $n=0$ case is trivial.
Assume the proposition
holds for $n-1$, i.e.,
\begin{align*}
\CT_{x}L(x_1,\ldots,x_n)D_{n-1}(x_1,\ldots,x_n;a_1,\ldots,a_n)=R(a_1,\ldots,a_n)\frac{d!}{a_1!\cdots
a_n!},
\end{align*}
where $d:=a_1+\cdots +a_n$, for any Laurent polynomial
$L(x_1,\ldots,x_n)$ independent of $a_i$. By linearity, it is
sufficient to show that \eqref{e-rationality} holds when
$L(x_0,x_1,\dots, x_n)$ is a monomial.  Define
\[f = f(a_0,a_1,\dots, a_n):=\CT_{x}x_0^{k_0} x_1^{k_1}\cdots x_{n}^{k_n}D_{n}(x_0,\ldots,x_n;a_0,a_1,\ldots,a_n).\]
We construct a rational function $R(a_0,a_1,\dots,a_n)$ so that
\begin{equation}
\label{e-s-rationality}
f(a_0,\dots,a_n)=R(a_0,a_1,\dots,a_n)\frac{(a_0+a_1+\cdots+a_n)!}{a_0!\,a_1!\cdots
a_n!}
\end{equation}
holds for all nonnegative
integers $a_0,a_1,\dots, a_n$.

First we show that for each
nonnegative integer $a_0$, there is a rational function $R_{a_0}(a_1,\dots,a_n)$ of $a_1,\dots,a_n$ such that
\begin{align}\label{ration3}
f(a_0,\dots, a_n)=R_{a_0}(a_1,\dots,a_n)\frac{d!}{a_1!\cdots
a_n!}.
\end{align}

By \eqref{expand1}, with
$i_l = a_0+a_l-j_l$ for $l=1,\dots,n$
we have
\begin{align*}
f=&\sum_{j_{1},\ldots,j_{n}}(-1)^{k_0} {a_0+a_{1}\choose j_{1}}
   \cdots{a_0+a_{n}\choose j_{n}}\cdot
   \CT_{x_{1},\ldots,x_{n}}\prod_{l=1}^{n}x_{l}^{a_0-j_{l}+k_l}
   D_{n-1}(x_1,\ldots,x_n;a_1,\ldots,a_n),
\end{align*}
where the sum ranges over all nonnegative integers $j_1,\dots,j_n$
such that $(a_0+a_1-j_1)+\cdots +(a_0+a_n-j_n)=d+k_0$, i.e.,
$j_1+\cdots+j_n=na_0-k_0.$

Therefore, by the induction hypothesis on $n$,
\begin{align*}
f&=\!\!\sum_{j_{1}+\cdots+j_{n}=na_0-k_0} (-1)^{k_0}{a_0+a_{1}\choose
j_{1}} \cdots{a_0+a_{n}\choose j_{n}}
R(j_1,\ldots,j_n;a_1,\ldots,a_n)\frac{d!}{a_1!\cdots a_n!},
\end{align*}
where for each $j_1,\dots,j_n$, $R(j_1,\ldots,j_n;a_1,\ldots,a_n)$
is a rational function of $a_1,\dots, a_n$ (which also depends on $a_0$ and $k_1,\dots, k_n$). Then \eqref{ration3}
holds with
\[R_{a_0}(a_1,\dots, a_n)=\sum_{j_{1}+\cdots+j_{n}=na_0-k_0} (-1)^{k_0}{a_0+a_{1}\choose j_{1}}
\cdots{a_0+a_{n}\choose j_{n}} R(j_1,\ldots,j_n;a_1,\ldots,a_n).\]

Now let $\beta_1=\beta_1(a_1,\dots, a_n),\dots,
\beta_r=\beta_r(a_1,\dots, a_n)$, where $r=(2^n-1)k$, be the linear
functions of $a_1,\dots, a_n$ of the form $\sigma(T)+j$, for
$T\subset[n]$ and $1\le j\le k$, where $k$ and $\sigma$ are as in
Lemma \ref{l-dys-root}.   By Lemma \ref{lemp}, $f(a_0,a_1,\dots, a_n)$, for fixed
$a_1,\dots, a_n$,  is a polynomial in $a_0$ of degree at most
$d+k_0$. Moreover, by Lemma \ref{l-dys-root},
 $f=0$ for $-a_0\in [d]\setminus \{\beta_1,\dots,\beta_r\}$.
 Thus there is a polynomial $p(a_0)$ of degree at most
$r+k_0$ (depending on $a_1,\dots, a_n$) such that
\begin{align}\label{ration2}
f(a_0,a_1,\dots, a_n)=\frac{(a_0+1)(a_0+2)\cdots (a_0+d)}{(a_0+\beta_1)(a_0+\beta_2)\cdots
(a_0+\beta_r)} p(a_0),
\end{align}
since  $f$ vanishes at the zeroes of the numerator factors that are not canceled by denominator factors.
Comparing with \eqref{ration3}, we obtain
\begin{equation*}
f=\frac{(a_0+d)!}{a_0!\,(a_0+\beta_1)(a_0+\beta_2)\cdots(a_0+\beta_r)}p(a_0)
=R_{a_0}(a_1,\ldots,a_n)\frac{d!}{a_1!\cdots a_n!}.
\end{equation*}
It follows that
\begin{align*}
p(a_0)=\frac{(a_0+\beta_1)\cdots(a_0+\beta_r)\,a_0!}{a_1!\cdots
a_n!\,(d+1)\cdots (d+a_0)}R_{a_0}(a_1,\ldots,a_n)= \frac{1}{a_1!\cdots
a_n!}\overline{R}_{a_0}(a_1,\ldots,a_n)
\end{align*}
for some rational function $\overline{R}_{a_0}(a_1,\ldots,a_n)$ of
$a_1,\dots,a_n$.

Applying  the Lagrange interpolation formula, we obtain that
\begin{align*}
p(a_0)=\sum_{l=0}^{r+k_0}p(l)\prod_{i=0,i\ne l}^{r+k_0}\frac{a_0-i}{l-i}
=\sum_{l=0}^{r+k_0}\frac{1}{a_1!\cdots
a_n!}\overline{R}_{l}(a_1,\ldots,a_n)\prod_{i=0,i\ne
l}^{r+k_0}\frac{a_0-i}{l-i}.
\end{align*}

So by \eqref{ration2}  we get
\begin{align*}
f&=\frac{(a_0+1)(a_0+2)\cdots
(a_0+d)}{(a_0+\beta_1)(a_0+\beta_2)\cdots (a_0+\beta_r)}
\sum_{l=0}^{r+k_0}\frac{1}{a_1!\cdots a_n!}
\overline{R}_{l}(a_1,\ldots,a_n)\prod_{i=0,i\ne l}^{r+k_0}\frac{a_0-i}{l-i}\\
&=R(a_0,\dots,a_n)\frac{(a_0+a_1+\cdots +a_n)!}{a_0!\,a_1!\cdots
a_n!},
\end{align*}
where
\[
R(a_0,\dots,a_n)=\frac{1}{(a_0+\beta_1)\cdots
(a_0+\beta_r)}\sum_{l=0}^{r+k_0}
\overline{R}_{l}(a_1,\ldots,a_n)\prod_{i=0,i\ne
l}^{r+k_0}\frac{a_0-i}{l-i}
\]
is a rational function of $a_0,\dots, a_n$.
 This completes the
induction.
\end{proof}

We note that the proof of Proposition \ref{p-rat} shows that the
denominator of $R(a_0,\dots, a_n)$ is a product of linear
polynomials of the form $a_{i_1}+\cdots+ a_{i_m}+j$, where $j$ is a
positive integer. This is consistent with the explicit formulas of
\cite{sills-zeilberger}, \cite{sills}, and \cite{LvXinZhou}.

\section{The Morris constant term identity \label{s-morris}}
 The proof of \eqref{e-morris} is
similar to that of \eqref{e-dyson}, so we omit some of the details.
We denote by $M'_n(a,b,k)$ the left side of \eqref{e-morris}.
\begin{lem}\label{l-allk}
For fixed $a\in \ZZ$ and $b\in \NN$, if $M'_{n}(a,b,k)=M_{n}(a,b,k)$
for $k\geq b$, then $M'_{n}(a,b,k)=M_{n}(a,b,k)$ for all $k\in
\mathbb{N}$.
\end{lem}
\begin{proof}
For fixed $a$ and $b$, by taking the constant term in $x_0$, we can
write $M'_{n}(a,b,k)$ as $\CT_{x} L D_{n-1}(x_1,\dots,x_n;
k,k,\dots,k)$, where $D_{n-1}$ and $L$ are as in Section \ref{s2}.
By Proposition~\ref{p-rat}, $M'_n(a,b,k)/M'_n(0,0,k)$ is a rational
function of $k$. It is straightforward to check that
$M_n(a,b,k)/M_n(0,0,k)$ is also rational in $k$. Note that
$M'_n(0,0,k)=M_n(0,0,k)$ follows from the equal parameter case of
the Dyson conjecture. Therefore, the hypothesis implies that
$M'_n(a,b,k)/M'_n(0,0,k)=M_n(a,b,k)/M_n(0,0,k) $ for all $k$. The
lemma then follows.
\end{proof}

\begin{proof}[Proof of Theorem \ref{t-morris}]
By setting $a_0=a, a_i=b$ for $i=1,\dots, n$ in Lemma \ref{lemp}, we
see that $M'_n(a,b,k)$ is a polynomial in $a$ of degree at most $bn$
for fixed $b$ and $k$ in $\mathbb{N}$. To see that $M_n(a,b,k) $
also has this property, we rewrite \eqref{mn} as
\begin{equation}\label{mn1}
M_n(a,b,k)=\prod_{l=0}^{n-1}\frac{(a+kl+1)(a+kl+2)\cdots(a+kl+b)(k(l+1))!}{(b+kl)!\,k!}.
\end{equation}
Moreover, it is easily seen that $M_n(a,b,k)$ vanishes if $a$ equals
one of the following values:
\begin{equation}\label{root}
\begin{array}{cccc}
  -1, & -2, & \ldots, & -b; \\
  -(k+1), & -(k+2), &\ldots, & -(k+b); \\
  \vdots & \vdots& \vdots & \vdots \\
  -[(n-1)k+1], & -[(n-1)k+2], & \ldots, & -[(n-1)k+b].
\end{array}
\end{equation}
Note that these values are distinct if $k\ge b$.

The theorem will follow from properties of polynomials, as in the
proof of Dyson's conjecture, if we can show that for $bn+1$ distinct
values of $a$, $M_n'(a,b,k)=M_n(a,b,k)$.

Lemma \ref{l-allk} reduces the problem to showing that
$M'_n(a,b,k)=M_n(a,b,k)$ if $k\ge b$. First we show that the
equality holds for $a=0$: From \eqref{hn}, we have
\begin{equation}\label{hn1}
\CT_{x_0}H(x_0,x_1,\ldots,x_n;0,b,k) = \prod_{1\leq i\neq j\leq n}
  \left( 1-\frac{x_i}{x_j} \right)^{\!\!k}.
\end{equation}
Thus by the equal parameter case of Dyson's conjecture,
$M'_n(0,b,k)={(nk)!}/{(k!)^n}$, which is equal to $M_n(0,b,k)$.

The remaining values are obtained from the following lemma, which completes the proof of Theorem
\ref{t-morris}.
\end{proof}

\begin{lem}\label{lem3.2}
For fixed nonnegative integers $b$ and $k\ge b$, $ \CT_x
H(x_0,x_1,\ldots,x_n;a,b,k)$ vanishes when $a$ equals one of the
values in \eqref{root}.
\end{lem}
\begin{proof}
We prove the contrapositive: Suppose that $h\in [nk]$ but the
constant term of
$$H(x_0,x_1,\dots,x_n;-h,b,k)= (-1)^{k\binom{n}{2}+nb}\frac{\prod_{1\le i<j\le n} (x_i-x_j)^{2k}}{\prod_{l=1}^n x_l^{(n-1)k+b}  \prod_{l=1}^n (1-x_l)^{h-b} } $$
 is not $0$. We shall show that
$(t-1)k+b< h \le tk$ for some $t$, i.e., $-h$ is not in \eqref{root}.

We apply Lemma \ref{l-ineq} with $u_{ij}=2k$, $m_i=(n-1)k+b$, and
$v_i = h - b$. Then  for some subset $T\subseteq[n]$ we have
\begin{align}
\sum_{\substack{i,j\in T\\ i<j}}2k &\le  \sum_{i\in T}(h-b) -|T| \label{e-morris-2.8d}\\[-20pt]
\intertext{and}\notag\\[-20pt]
\sum_{1\le i<j\le n}\!\!2k &\ge \sum_{\substack{i,j\in T\\
i<j}}2k
 +\sum_{i\in \T}\bigl((n-1)k+h\bigr).
\label{e-morris-2.9d}
\end{align}
Let $t=|T|$. Then \eqref{e-morris-2.8d} may be written as
\begin{equation*}
(t-1)tk\le t(h-b)-t,
\end{equation*}
and this implies that for $t\ne0$,
\begin{equation}
\label{morris4} (t-1)k+b< h.
\end{equation}
But \eqref{morris4} also holds for $t=0$, since $h\ge
1>b-k$.

Similarly, \eqref{e-morris-2.9d} gives
\begin{equation*}
(n-1)nk \ge t(t-1)k+(n-t)((n-1)k+h),\end{equation*}
which simplifies to  $t(n-t)k\ge(n-t)h$,
so for $t\ne n$,
\begin{equation}
\label{morris5a} h\le tk.
\end{equation}
But \eqref{morris5a} also holds for $t=n$, since $h\le nk$. Thus by \eqref{morris4} and \eqref{morris5a},
\[ (t-1)k+b<h\le tk,\]
which completes the proof.
\end{proof}

We note that it is possible to  handle the case $k<b$ directly
without applying Proposition~\ref{p-rat} by using the following
fact: $z_0$ is a root of a polynomial $P(z)$ with multiplicity $r$
if and only if $\displaystyle\frac{d^i}{dz^i}P(z_0)=0$ for
$i=0,1,\dots, r-1$. For instance, we can find roots of $M'_n(a,b,k)$
with multiplicity at least $2$ by considering the constant term of
\begin{equation*}
\frac{\partial}{\partial a} H(x_0,x_1,\ldots,x_n;a,b,k)=
\sum_{s=1}^n \ln\left(1-\frac{x_s}{x_0}\right)\prod_{l=1}^n
\left(1-\frac{x_l}{x_0}\right)^{\!\!a}
\left(1-\frac{x_0}{x_l}\right)^{\!\!b}\prod_{1\leq i\neq j\leq n}
  \left( 1-\frac{x_i}{x_j} \right)^{\!\!k}.
\end{equation*}

\section{The Aomoto constant term identity\label{s-aomoto}}

In this section we will prove Aomoto's identity using our elementary
approach. First we note that if $m\leq 0$ or $m\geq n$, then
\eqref{extmorris} reduces to the Morris identity, 
so we assume here
that  $1\le m \le n-1$. The proof is similar to that of the Morris
identity but is more complicated. We provide only the details of the
key points.

In contrast with the Morris identity, it is not easy to show that
\eqref{extmorris} holds when $a=0$. So instead of proving equality
at a $bn+1$st point, we show that both sides of \eqref{extmorris}
have the same leading coefficients as polynomials in $a$.

\begin{prop}\label{exmprop}\ \\[-17pt]
\begin{itemize}
 \item[\bf{1.}] Both sides of \eqref{extmorris} are polynomials in $a$ of degree at most $bn$.
\item[ \bf{2.}] The left side and the right side of \eqref{extmorris} have the same leading coefficients in $a$.
\item[\bf{3.}] The right side of \eqref{extmorris} vanishes when $a$ equals one of the values in
the following table.
\begin{equation}\label{emroot}
\begin{array}{cccc}
  -1, & -2, & \ldots, & -b; \\
  -(k+1), & -(k+2), &\ldots, & -(k+b); \\
  \vdots & \vdots& \vdots & \vdots \\
  -[(n-m-1)k+1], & -[(n-m-1)k+2], & \ldots, & -[(n-m-1)k+b];
  \\[-2pt]  
 \multicolumn{4}{c}{\quad\hrulefill\quad}\\
 \\[-10pt]
-[(n-m)k+2], & -[(n-m)k+3], & \ldots, & -[(n-m)k+b+1]; \\
 \vdots & \vdots& \vdots & \vdots \\
-[(n-1)k+2], & -[(n-1)k+3], & \ldots, & -[(n-1)k+b+1]. \\
\end{array}
\end{equation}
\end{itemize}
\end{prop}

\begin{proof}[Proof of Proposition {\rm\ref{exmprop}} \rm{(sketch)}]
As with the Morris identity, parts 1 and 3 are straightforward. To
show part 2, we rewrite the right side of \eqref{extmorris} as
\[
\prod_{l=0}^{n-1}\frac{[a+kl+\chi(l\geq
n-m)+1]\cdots[a+kl+\chi(l\geq n-m)+b](k(l+1))!}{(b+kl)!\,k!},
\]
whose leading coefficient is now clearly
\begin{equation}\label{emlpl}
\prod_{l=0}^{n-1}\frac{(k(l+1))!}{(b+kl)!\,k!}.
\end{equation}

On the other hand, a calculation similar to that in Lemma \ref{lemp}
shows that the leading coefficient of the left side of
\eqref{extmorris} equals
\begin{align*}
\frac{1}{(nb)!}\CT_x(x_1+\cdots+x_n)^{nb}\prod_{l=1}^{n}x_{l}^{-b}\prod_{1\le
i\ne j \le n} \left(1-\frac{x_i}{x_j}\right)^{\!\!k},
\end{align*}
which is equal to \eqref{emlpl} by Corollary \ref{c-leadingM}.
\end{proof}

As in the proof of the Morris identity, we may assume $k$ to be
sufficiently large by Proposition \ref{p-rat}. Then we can complete
the proof of the Aomoto identity by the following lemma.
\begin{lem}\label{lem-ext-3.2}
For fixed nonnegative integers $b$, $k\ge b$, and $m \in [n]$, if
$a$ equals one of the values in \eqref{emroot}, then $ \CT_x
A_m(x_0,x_1,\ldots,x_n;a,b,k)$ vanishes.
\end{lem}
\begin{proof}
We prove the contrapositive: Suppose that $h\in [nk+1]$ but the
constant term of
$$A_m(x_0,x_1,\dots,x_n;-h,b,k)= (-1)^{k\binom{n}{2}+nb}\frac{\prod_{1\le i<j\le n} (x_i-x_j)^{2k}}{\prod_{l=1}^n x_l^{(n-1)k+b}  \prod_{l=1}^n (1-x_l)^{h-b-\chi(l\le m)} } $$
 is not equal to $0$. We shall show that
$(t-1)k+b+1\le  h \le tk$ for some $t$ with $1\le t< n-m$, or
$(t-1)k+b+1\le  h \le tk+1$ for $t=n-m$, or $(t-1)k+b+2\le h \le
tk+1$ for some $t$ with $n-m\le t\le n$. That is, $-h$ is not in
\eqref{emroot}.

We apply Lemma \ref{l-ineq} with $u_{ij}=2k$, $m_i=(n-1)k+b$, and
$v_i = h - b-\chi(i\le m)$. Then  for some subset $T\subseteq[n]$ we
have
\begin{align}
\sum_{\substack{i,j\in T\\ i<j}}2k &\le  \sum_{i\in T}(h-b-\chi(i\le m))-|T| \label{e-extmorris-2.8d}\\[-20pt]
\intertext{and}\notag\\[-20pt]
\sum_{1\le i<j\le n}\!\!2k &\ge \sum_{\substack{i,j\in T\\
i<j}}2k
 +\sum_{i\in \T}\bigl((n-1)k+h -\chi(i\le m) \bigr).
\label{e-extmorris-2.9d}
\end{align}
Let $t=|T|$. Then \eqref{e-extmorris-2.8d} may be written as
\begin{equation*}
(t-1)tk\le t(h-b)-t-\sum_{i\in T} \chi(i\le m),
\end{equation*}
and this implies that for $t\ne0$,
\begin{equation}
\label{extmorris4} (t-1)k+b+
2-\chi(T \cap [m]=\varnothing )\le h.
\end{equation}
But \eqref{extmorris4} also holds for $t=0$, since $h\ge
1\ge b-k+1$.

Similarly, \eqref{e-extmorris-2.9d} gives
\begin{equation*}
(n-1)nk \ge t(t-1)k+(n-t)((n-1)k+h)-\sum_{i\in \T} \chi(i\le m).
\end{equation*}
Taking all terms in the $k$ to the left gives
$$t(n-t)k \ge (n-t)h -\sum_{i\in \T} \chi(i\le m),
  $$
so for $t\ne n$,
\begin{equation}
\label{extmorris5a} h\le tk
+\chi(\T \subseteq [m]).
\end{equation}
But \eqref{extmorris5a} also holds for $t=n$, since $h\le nk+1$.
Thus by \eqref{extmorris4} and \eqref{extmorris5a},
\[ (t-1)k+b+2-\chi(T \cap [m]=\varnothing )
\le h\le tk+\chi(\T \subseteq [m]).
\]
Now according to the three cases $t<n-m$, $t=n-m$, and $t>n-m$, the
minimum values of $-\chi(T\cap [m]=\varnothing)$ are $-1$, $-1$, and $0$,
respectively, and  the maximum values of $\chi(\T \subseteq [m])$ are $0,$ $1$,
and $1$, respectively. This completes the proof.
\end{proof}

\section{On the Forrester conjecture}

We can apply our method to Forrester's constant term to obtain some
partial results. It is routine to obtain the following.
\begin{prop}\label{pro}\ \\[-17pt]
\begin{itemize}
 \item[\bf{1.}] Both sides of \eqref{e-forrester} are polynomials in $a$ of degree at
most $bn$.
\item[ \bf{2.}] If $a=0$, then the left side of \eqref{e-forrester} is equal to the right side of \eqref{e-forrester}.
\item[\bf{3.}] The right side of \eqref{e-forrester} vanishes when
$a$ equals one of the values in the following table.
\begin{equation}\label{froot}
\begin{array}{cccc}
  -1, & -2, & \ldots, & -b; \\
  -(k+1), & -(k+2), &\ldots, & -(k+b); \\
  \vdots & \vdots& \vdots & \vdots \\
  -[(n_0-1)k+1], & -[(n_0-1)k+2], & \ldots, & -[(n_0-1)k+b]; \\[-2pt]
 \multicolumn{4}{c}{\quad\hrulefill\quad}\\
 \\[-10pt]
-(n_0k+1), & -(n_0k+2), & \ldots, & -(n_0k+b); \\
-[(n_0+1)k+2], & -[(n_0+1)k+3], & \ldots, & -[(n_0+1)k+b+1]; \\
 \vdots & \vdots& \vdots & \vdots \\
-[(n-1)k+n_1], & -[(n-1)k+n_1+1], & \ldots, & -[(n-1)k+n_1+b-1].
\end{array}
\end{equation}
\end{itemize}
\end{prop}
Note that the values in \eqref{froot} are distinct if $k\geq b.$

Therefore, by applying Proposition \ref{p-rat}, Forrester's
conjecture would be established if we could show that for
sufficiently large $k$, the left side of \eqref{e-forrester}
vanishes when $a$ equals any value in \eqref{froot}. However, we are
only able to show that it vanishes for some of these values. Denote by
$F_{n_0}(x;a,b,k)$ the left side of \eqref{e-forrester}. We obtain
the following.
\begin{lem}\label{p-roots}
Assume $k$ is sufficiently large. For $t$ with $0\le t\le n-1$, let
$M:=\min\{n_{1},t\}$. If $a=-h$ with $h$ satisfying the conditions
\begin{alignat}{2}
tk+C_{1}+1&\leq h\leq tk+b, & \quad&\text{if $0\leq t\leq n_{0}$,} \label{msmall}\\
tk+C_{2}+1&\leq h \leq  tk+b+C_{3}, &\quad&\text{if $n_{0}+1\leq
t\leq n-1$,}\label{mlarge}
\end{alignat}
where\\[-20pt]
\begin{align*}
C_{1}&=\begin{cases}
    \left\lfloor\frac{n_{1}^2}{4(n-t)}\right\rfloor,& \text{if $\frac{n_{1}}{2}\leq M$,}\\[7pt]
    \left\lfloor\frac{M(n_{1}-M)}{n-t}\right\rfloor,&\text{if $\frac{n_{1}}{2}>M$,}
\end{cases} \\
C_{2}&=\begin{cases}
    \left\lfloor\frac{M(n_{1}-M)}{n-t}\right\rfloor,&\text{if $\frac{n_{1}}{2}\geq M$,}\\[7pt]
    \left\lfloor\frac{n_{1}^2}{4(n-t)}\right\rfloor,&\text{if  $t-n_{0}< \frac{n_{1}}{2}< M$,}\\[2pt]
     t-n_{0}, &\text{if $\frac{n_{1}}{2}\leq t-n_{0}$,}
\end{cases} \\
C_{3} &=\left\lceil \frac{(t-n_{0}+1)(t-n_{0})}{t+1}\right\rceil,
\end{align*}
then $F_{n_0}(x;a,b,k)$ vanishes.
\end{lem}
\begin{proof}
We prove the contrapositive: Suppose that $h\in [nk+n_1]$ but the
constant term of
$$ F_{n_0}(x;-h,b,k)=\pm \frac{\prod_{1\le i<j\le n} (x_i-x_j)^{2(k+\chi_{ij}^{n_0})}}{\prod_{l=1}^n x_l^{(n-1)k+b+(n_1-1)\chi(l>n_0)}  \prod_{l=1}^n (1-x_l)^{h-b} }, $$
where $\chi_{ij}^{n_0}=\chi(i>n_0)\chi(j>n_0)$, is not equal to $0$.
We shall obtain conditions on $h$ from which the lemma follows.

We apply Lemma \ref{l-ineq} with $u_{ij}=2(k+\chi_{ij}^{n_0})$,
$m_i=(n-1)k+b+(n_1-1)\chi(i>n_0)$, and $v_i = h - b$. Then  for some
subset $T\subseteq[n]$ we have
\begin{align}
\sum_{\substack{i,j\in T\\ i<j}}(2k+2\chi_{ij}^{n_0}) &\le  \sum_{i\in T}(h-b) -|T| \label{e-forrester-2.8d}\\[-20pt]
\intertext{and}\notag\\[-20pt]
\sum_{1\le i<j\le n}\!\!(2k+2\chi_{ij}^{n_0}) &\ge \sum_{\substack{i,j\in T\\
i<j}}(2k+2\chi_{ij}^{n_0})
 +\sum_{i\in \T}\bigl((n-1)k+h+(n_1-1)\chi(i>n_0)  \bigr).
\label{e-forrester-2.9d}
\end{align}
Let $t=|T|$ and assume that in $T$ there are $t_0$ elements less than or equal to
$n_0$ and $t_1$ elements greater than $n_0$. Then \eqref{e-forrester-2.8d} may be
written as
\begin{equation*}
(t-1)tk +t_1(t_1-1) \le t(h-b)-t,
\end{equation*}
where $\max\{0,t-n_0\}\leq t_1\leq \min\{n_1,t\}$ by its definition,
and the above equation implies that for $T\ne \varnothing$,
\begin{equation}\label{forrester4}
(t-1)k+b+1+\frac{t_1(t_1-1)}{t}\leq  h.
\end{equation}
But \eqref{forrester4} also holds for $T=\varnothing$ if
$t_1(t_1-1)/t$ is taken as $-1$ when $t=0$ (hence $ t_1=0$), since
$h\ge 1>b-k$.

Similarly, \eqref{e-forrester-2.9d} gives
\begin{align*}
&(n-1)n k+(n_1-1)n_1 \ge
t(t-1)k+t_1(t_1-1)+(n-t)\big((n-1)k+h\big)+(n_1-1)(n_1-t_1).
\end{align*}
Taking all terms in the $k$ to the left gives
$$t(n-t)k\ge (n-t)h-t_1(n_1-t_1), $$
 so for $T\ne [n]$,
\begin{equation}
\label{forrester5a} h\le tk+\frac{t_1(n_1-t_1)}{n-t}.
\end{equation}
But \eqref{forrester5a} also holds for $T=[n]$ if
$t_1(n_1-t_1)/(n-t)$ is taken as $n_1$ when $t=n$ (hence $t_1=n_1$),
since $h\le nk+n_1$.

Thus by \eqref{forrester4} and
\eqref{forrester5a},
\[h \in I(t,t_1):=\left[(t-1)k+b+1+\frac{t_1(t_1-1)}{t}, tk+\frac{t_1(n_1-t_1)}{n-t}\right].\]
It follows that $F_{n_0}(x;-h,b,k)$ vanishes if
$$h\in [nk+n_1] \setminus \bigcup_{ t_1\le t} I(t,t_1),  $$
where $t$ ranges from $0$ to $n$ and $t_1$ ranges from
$\max\{0,t-n_0\}$ to $\min\{n_1,t\}$.

The above condition can be simplified further: For $0\leq t\leq n-1$
if
\begin{equation}\label{forrester7}
tk+\Big\lfloor\frac{r_1(n_1-r_1)}{n-t}\Big\rfloor+1\leq h\leq
tk+b+\Big\lceil\frac{r_2(r_2-1)}{t+1}\Big\rceil
\end{equation}
holds for every $r_1$ with $\max\{0,t-n_0\}\leq r_1\leq
\min\{n_1,t\}$ and $r_2$ with $\max\{0,t-n_0+1\}\leq r_2\leq
\min\{n_1,t+1\}$, then $F_{n_0}(x;-h,b,k)$ vanishes. This is because
when $k$ is sufficiently large, the  left and right endpoints of
$I(t,r_1)$ are always to the left of the corresponding endpoints of $I(t+1,r_2)$ for any
$r_1$ and $r_2$ in their range. Therefore, after removing the
intervals $\bigcup_{r_1} I(t,r_1)$ and $\bigcup_{r_2} I(t+1,r_2)$,
each remaining value of $h$  belongs to an open interval
(possibly empty), from the right endpoint of $I(t,r_1)$ to the left
endpoint of $I(t+1,r_2)$ for some $r_1$ and $r_2$.

By analyzing the extreme values of \eqref{forrester7} among the
range of $r_1$ and $r_2$, it is straightforward to obtain
\eqref{msmall} and \eqref{mlarge}.
\end{proof}

\begin{rem}\label{rem1}
The first two lines of {\rm\eqref{froot}} are always implied by
Lemma $\ref{p-roots}$ to be roots for $n_0\geq 1.$ This follows
easily by checking the cases $t=0$ and $t=1$.
\end{rem}

\begin{cor}\label{cor1}
Conjecture $\ref{conjecture}$ holds in the extreme cases $n_{1}=2$
and $n_{1}=n-1$.
\end{cor}
\begin{proof}
We verify this directly by Lemma \ref{p-roots}.

If $n_{1}=2$, then $n_{0}=n-2$. The first two lines of \eqref{froot}
are roots by Remark \ref{rem1}. If $2\leq t\leq n_{0}=n-2$, then
$M=\min\{2,t\}=2$, and $C_{1}=0$. Thus we obtain the range $tk+1\leq
h\leq tk+b$, which is consistent with the $(t+1)$st line of
\eqref{froot}. If $t=n-1$, then $M=2$, and $C_{2}=C_{3}=1$. Thus we
obtain the range $(n-1)k+2\leq h\leq (n-1)k+b+1$, which is
consistent with the $n$th line of \eqref{froot}. Therefore, Lemma
\ref{p-roots} implies that all values of $a$ in \eqref{froot} are
roots, and Forrester's conjecture holds in this case.

If $n_{1}=n-1$, then $n_{0}=1$. The cases $t\leq n_{0}=1$ of
\eqref{froot} are dealt with in Remark \ref{rem1}. If $2\leq t \leq
n-1$, then $M=t$,
$C_{3}=\big\lceil\frac{t(t-1)}{t+1}\big\rceil=\big\lceil\frac{(t+1)(t-1)-(t-1)}{t+1}\big\rceil=t-1$,
and all three cases of $C_{2}$ are equal to  $t-1$. Thus we
obtain the range $tk+t\leq h\leq tk+b+t-1$, which is consistent with
that of \eqref{froot}. As in the case $n_1=2$, Forrester's
conjecture holds.
\end{proof}

\begin{prop}\label{p-n5}
Forrester's conjecture holds when $n\leq 5$.
\end{prop}
\begin{proof}
The cases $n\leq 4$ are consequences of Corollary \ref{cor1}; the
case $n=5$ can be verified by Lemma \ref{p-roots}.
\end{proof}

Further routine calculations by Lemma \ref{p-roots} gives us the
following table:
\begin{align}\label{diagram1}
\begin{array}{ccccccccc}
n_{1}= &2 &3 &4 &5 &\cdots &n-3    &n-2 &n-1\\
M_{r}=     &0  &1 &4 &8 &\cdots &2n_{1}+2\lfloor
\frac{n_{1}}{2}\rfloor-5 &n_{1}+\lfloor \frac{n_{1}}{2}\rfloor-6 &0
\end{array}
\end{align}
where $M_{r}$ is an upper bound for the number of \emph{missing
roots} in \eqref{froot}, i.e., roots that are not implied by Lemma
\ref{p-roots}. For brevity, we verify  in detail only the case $n_1=3$;
the other cases are similar.

For $n_{1}=3$, we have $n_{0}=n-3$. The cases $t=0,1$ are guaranteed
by Remark \ref{rem1}, so henceforth we will always assume $t\ge 2$
and (by Proposition \ref{p-n5}) $n\ge 6$. If $t=2$, then $M=2$ and
$\frac{n_{1}}{2}=3/2<M$. Therefore
$C_{1}=\Big\lfloor\frac{3^2}{4(n-2)}\Big\rfloor=0$, for $n\geq 5$.
If $3\leq t\leq n-3$, then $M=3$, and
$C_{1}=\Big\lfloor\frac{n_{1}^2}{4(n-t)}
\Big\rfloor=\Big\lfloor\frac{9}{4(n-t)}\Big\rfloor=0$. If $t=n-2$,
then $M=3$ and $t-n_{0}=1$. It follows that $t-n_{0}<
\frac{n_{1}}{2}< M$, $C_{2}=\Big \lfloor \frac{n_{1}^2}{4(n-t)}\Big
\rfloor=1$, $C_{3}=1$. If $t=n-1$, then $M=3$. This implies that
$\frac{n_{1}}{2}\leq t-n_{0}$, $C_{2}=t-n_{0}=2$, and $C_{3}=1$.

In conclusion, in the case $n_{1}=3$, only one root $-[(n-1)k+b+2]$
is not implied by Lemma \ref{p-roots}.

\begin{cor}
Conjecture $\ref{conjecture}$ holds in the case $n_{1}=3$.
\end{cor}

\begin{proof}

As we have just seen, for $n_1=3$ we are missing only one root. But by \cite{baker-forr}, we know that the
$q$-generalization of \eqref{e-forrester} holds when $a=k$. Thus if
we let $q\rightarrow 1$ in this result, we get $a=k$ as our $(bn+1)$st
point.
\end{proof}

We conclude this paper by the following observation. Let us take the
Forrester constant term as an example. In the proof of Lemma
\ref{l-ineq}, we made the expansion
\[\prod_{1\le i< j \le n }(y_i-y_j)^{2k}\!\! \prod_{n_0+1\le i<j\le n}
(y_i-y_j)^2=\sum_{\gamma } c_\gamma \, y_1^{\gamma_1} \cdots
y_n^{\gamma_n}, \] where $y_i=1-x_i$, and try to show that the
constant term associated with $y_1^{\gamma_1} \cdots y_n^{\gamma_n}$
is equal to $0$ for each $\gamma$. However, it would be sufficient
to show that for each $\gamma$, either the associated constant term
is $0$ or $c_\gamma=0$ (after cancellation). We \emph{conjecture}
that in our approach to Forrester's conjecture, $c_\gamma=0$ when
Lemma \ref{l-ineq} does not apply. We have checked our conjecture
for $n\le 6$ and $k\le 3$.

\vspace{.2cm} \noindent{\bf Acknowledgments.} The authors would like
to thank the referees for helpful suggestions to improve the
presentation. Lun Lv and Yue Zhou would like to acknowledge the
helpful guidance of their supervisor William Y.C. Chen. The third
author was supported by the 973 Project, the PCSIRT project of the
Ministry of Education, the Ministry of Science and Technology and
the National Science Foundation of China.

\section{Appendix: Consequences of the polynomial approach \label{appendix}}

From  Lemma \ref{lemp} and its proof, we can deduce the following result.

\begin{cor}\label{corleading}
Let $d$ and $Q$ be as in Lemma \ref{lemp} with $k_0=0$. Then the
leading coefficient of $Q(a_0,\dots,a_n)$ in $a_0$ is
\begin{align}
\frac{1}{d!}\CT_{x}(x_{1}+\cdots+x_{n})^{d}\prod_{l=1}^{n}x_{l}^{-a_{l}}
L(x_1,\dots,x_n),\label{e-leading}
\end{align}
and the second leading coefficient of $Q(a_0,\dots,a_n)$ in $a_0$ is
\begin{align}
\CT_{x}\left(\sum_{l=1}^{n}\frac{a_{l}x_{l}}{(d-1)!}\Big(\sum_{i=1}^n
x_{i}\Big)^{\!d-1}-
\frac{1}{2}\sum_{l=1}^{n}\frac{x_{l}^2}{(d-2)!}\Big(\sum_{i=1}^n
x_{i}\Big)^{\!d-2}\right)
\prod_{l=1}^{n}x_{l}^{-a_{l}}L(x_{1},\ldots,x_{n}).\label{e-sleading}
\end{align}
\end{cor}
\begin{proof}
Taking the leading coefficient of \eqref{expand1} in $a_0$ gives
\begin{align*}
\sum_{i_{1},\ldots,i_{n}}\frac{1}{i_{1}!\,i_{2}!\cdots
i_{n}!}\CT_{x}\prod_{l=1}^{n} \kern -5pt&\kern5pt
x_{l}^{i_{l}-a_{l}}L(x_{1},\ldots,x_{n})\nonumber \\
=&\frac{1}{d!}\CT_{x}\sum_{i_{1},\ldots,i_{n}}\frac{d!\, x_{1}^{i_{1}}x_{2}^{i_{2}}\cdots x_{n}^{i_{n}}}{i_{1}!\, i_{2}!\cdots i_{n}!}\prod_{l=1}^{n}x_{l}^{-a_{l}}L(x_{1},\ldots,x_{n})\nonumber\\
=&\frac{1}{d!}\CT_{x}(x_{1}+\cdots+x_{n})^{d}\prod_{l=1}^{n}x_{l}^{-a_{l}}
L(x_1,\dots,x_n).
\end{align*}

Taking the second leading coefficient of \eqref{expand1} gives
\begin{align*}
&\CT_{x}\sum_{i_{1},\ldots,i_{n}}\sum_{l=1}^{n}\frac{i_{l}(2a_{l}-i_{l}+1)}{2}
\frac{x_{1}^{i_{1}}x_{2}^{i_{2}}\cdots
x_{n}^{i_{n}}}{i_{1}!\,i_{2}!\cdots i_{n}!}
\prod_{l=1}^{n}x_{l}^{-a_{l}}L(x_{1},\ldots,x_{n}),
\end{align*}
which can be rewritten as \eqref{e-sleading}.
\end{proof}

Applying Corollary \ref{corleading} to the Dyson conjecture gives
the following identity, which appeared in \cite[Corollary~5.4]{xin}.
\begin{cor}
\begin{align}\label{D-leading}
\CT_{x}(x_{1}+\cdots+x_{n})^{a_1+\cdots+a_n}\prod_{l=1}^{n}x_{l}^{-a_{l}}\prod_{1\leq
i\neq j\leq
n}\left(1-\frac{x_{i}}{x_{j}}\right)^{\!\!a_{j}}=\frac{(a_1+\cdots+a_n)!}{a_{1}!\,a_{2}!\cdots
a_{n}!}.
\end{align}
\end{cor}
We omit the formula for the second leading coefficient, which is more complicated.

Applying Corollary \ref{corleading} to  Morris's identity gives the
following result, which is needed for the proof of Aomoto's
identity.  We remark that
\eqref{M-leading} was also obtained in
\cite[Proposition~2.2]{hou-sun} through a complicated calculation.

\begin{cor}\label{c-leadingM}
\begin{align}\label{M-leading}
\CT_{x}(x_{1}+\cdots+x_{n})^{nb}\prod_{l=1}^{n}x_{l}^{-b}\prod_{1\leq
i\neq j\leq
n}\left(1-\frac{x_{i}}{x_{j}}\right)^{\!\!k}=(nb)!\prod_{l=0}^{n-1}\frac{(k(l+1))!}{(b+kl)!\,k!}.
\end{align}
\small{\begin{align}\label{M-Sleading} \CT_{x}
\frac{x_{1}^2\Big(\sum_{i=1}^n x_{i}\Big)^{\!nb-2}}{
\prod_{l=1}^{n}x_{l}^{b}} \prod_{1\leq i\neq j\leq
n}\left(1-\frac{x_{i}}{x_{j}}\right)^{\!\!k}
&=-b\big(nk-k-b+1\big)(nb-2)!\prod_{l=0}^{n-1}\frac{(k(l+1))!\,
}{(b+kl)!\,k!}.
\end{align}}
\end{cor}

The proofs of \eqref{D-leading}, \eqref{M-leading} and
\eqref{M-Sleading} are straightforward.

\renewcommand\emph{}
\bibliographystyle{amsplain}

\providecommand{\bysame}{\leavevmode\hbox
to3em{\hrulefill}\thinspace}
\providecommand{\MR}{\relax\ifhmode\unskip\space\fi MR }
\providecommand{\MRhref}[2]{%
  \href{http://www.ams.org/mathscinet-getitem?mr=#1}{#2}
} \providecommand{\href}[2]{#2}

\end{document}